\documentclass[preprint,3p, review, 10pt]{elsarticle}
\usepackage{mathrsfs}
\usepackage{amssymb}
\usepackage{amsmath}
\usepackage{caption}
\usepackage{graphicx}
\usepackage{latexsym}
\usepackage{amsthm}
\usepackage{multirow}
\usepackage{lineno}
\usepackage[colorlinks, citecolor=red]{hyperref}
\usepackage{cleveref} 
\usepackage{booktabs}
\usepackage{multirow}
\usepackage{threeparttable}

\crefname{theorem}{Theorem}{Theorems}
\crefname{lemma}{Lemma}{Lemmas}
\journal{Mathematical Methods in the Applied Sciences}
\numberwithin{equation}{section}

\newtheorem{theorem}{Theorem}[section]
\newtheorem{lemma}{Lemma}[section]
\newtheorem{definition}{Definition}[section]

\newtheorem{remark}{Remark}[section]

\allowdisplaybreaks
\begin{document}

\begin{frontmatter}



\title{Global existence and blow-up for a fractional wave equation with logarithmic nonlinearity and strong damping} 


\author{Guangyang Fu} 
\ead{107552300585@stu.xju.edu.cn}
\author{Jiang Zhou\corref{cor1}} 
\ead{zhoujiang@xju.edu.cn}
\address{College of Mathematics and System Sciences, Xinjiang University, Urumqi 830046, PR China}

\cortext[cor1]{Corresponding author}
\begin{abstract}
This paper studies a fractional strongly damped wave equation with logarithmic nonlinearity. Based on local well‑posedness of weak solutions via the Galerkin approximation, we employ the potential well theory to derive the existence of global solutions under subcritical and critical initial energy conditions. Moreover, infinite time blow-up phenomena are obtained under subcritical, critical, and arbitrarily high supercritical initial energy conditions.
\end{abstract}


\end{frontmatter}



	\section{Introduction}

In this paper, we study the initial boundary value problem (IBVP) of a fractional wave equation with logarithmic nonlinearity and strong damping
\begin{eqnarray} \label{1.1}
	\begin{cases}
		u_{tt}+(-\Delta)^s u_t+(-\Delta)^s u= u \ln|u|, & x\in \Omega,\, t>0, \\
		u(x,t)=0, & x\in \mathbb{R}^n\setminus \Omega,\, t>0, \\
		u(x,0)=u_0(x),\quad u_t(x,0)=u_1(x), & x\in \Omega,
	\end{cases}
\end{eqnarray}
where $0<s<1$, $\Omega\subset \mathbb{R}^n$ ($n\in \mathbb{N}$ and $n > 2s$) is a bounded domain with smooth boundary $\partial \Omega$. The fractional Laplacian operator $(-\Delta)^s$ is defined by 
\begin{eqnarray} \label{b1.2}
	(-\Delta)^s u(x)=C(n,s) \lim_{\varepsilon \to 0^+} \int_{\mathbb{R}^n \setminus B_\varepsilon(x)} \frac{u(x)-u(y)}{|x-y|^{n+2s}} dy,\quad x\in \Omega, \quad (\text{see } \cite{2012DNP}),
\end{eqnarray}
where the normalization constant is given by
\[
C(n,s) = \left(\int_{\mathbb{R}^n} \frac{1-\cos \left(\zeta_1\right)}{|\zeta|^{n+2 s}} d \zeta\right)^{-1} = - \frac{ 4^s \Gamma(\frac{n}{2}+s)}{\pi^{n / 2} \Gamma(-s)}, \quad \text{(see \cite[Lemma~5.1]{2010ST})}.
\]

The fractional derivative is a nonlocal operator with integro-differential properties. Unlike its classical local counterpart, its value at an arbitrary point is determined by both the local function values near the point and the function values across the full spatial region. Specifically, spatial nonlocality arises from arbitrarily large spatial jumps, and temporal nonlocality comes from history-dependent effects in dynamics induced by anomalously extended waiting times \cite{2021ZW}. Currently, fractional partial differential equations have been extensively employed in a wide range of disciplines, including finance \cite{2023AG}, thin obstacle problems \cite{2007S}, crystal dislocation \cite{1997T}, anomalous diffusion \cite{2018IS}, flame propagation \cite{2012CMS}, and the ultra-relativistic limit in quantum mechanics \cite{1986FL}. We consider the fractional Laplacian wave equation
\begin{eqnarray} \label{1.2}
	u_{tt}+(-\Delta)^s u=f(u). 
\end{eqnarray}
For the homogeneous case $f(u)=0$, in 2012, Erbay et al. \cite{2012EEM} analyzed the Cauchy problem for a general nonlinear peridynamic equation in one dimension. In 2018, this problem was later generalized to $N$-dimensional spaces by Coclite et al. \cite{2018CDMV}. For the linear source term $f(u)=b(x, t)$, in 2014, Emmrich and Dimitri \cite{2014ED} considered the Cauchy problem in two and three spatial dimensions. For the polynomial source term $f(u)= \pm|u|^{p-1} u$, in 2017, Dinh \cite{2017D} addressed the local well-posedness of weak solutions to the Cauchy problem. In 2025, Lin and Xu \cite{2025LX} further generalized equation \eqref{1.2} to the variable-order fractional case with the operator $(-\Delta)^{s(\cdot)}$ and variable-exponent nonlinearity $f(u) = |u|^{p(\cdot)-2} u$, and studied the global well-posedness of solutions. In addition, Kirchhoff-type wave equations involving the fractional Laplacian have also attracted research attention in recent years \cite{2025CHX,2024LT,2025SX}.

Logarithmic nonlinearity appears in various areas of physics, including quantum mechanics and inflation cosmology \cite{2008BG,1975BBM,1976BBM,2003B}. It also arises in nuclear physics, optics, and geophysics \cite{1980CH,2003DM,2009G,2000K}. Due to these special physical meanings, research on evolution equations with logarithmic nonlinearity has attracted great attention. Readers are referred to 
\cite{2020CX,2020DS,2019DZ,2018HH,2019LASHX,2019LWX,2024LB,2025SX,2023W,2024WW,2022Y} and the references therein. In 2025, Ding and Zhou \cite{2025DZ} considered the following initial boundary value problem for a $p$-Laplacian type wave equation ($1<p \leq 2$)
\begin{eqnarray} \label{1.3}
			\begin{cases}u_{t t}-\Delta u_t-\Delta_p u=|u|^{r-2} u \log |u|, & x \in \Omega, t>0, \\ u(x, t)=0, & x \in \partial \Omega, t>0, \\ u(x, 0)=u_0(x), u_t(x, 0)=u_1(x), & x \in \Omega. \end{cases}
\end{eqnarray}
They first proved the local existence of weak solutions via the Galerkin method. Using potential well theory \cite{1968S} and differential inequality techniques, they established the global existence of solutions. They also proved finite time and infinite time blow-up results under subcritical, critical, and supercritical initial energy conditions. Furthermore, they derived upper and lower bounds for the blow-up time of finite time blow-up solutions.
In \cite{2025DZ}, the authors used the logarithmic Sobolev inequality as a key tool to establish the main results for the case \(p=r=2\) of the IBVP \eqref{1.3}. We recall this inequality as follows: for any $\phi \in H_0^1(\Omega)\setminus\{0\}$ and $a>0$, there holds
\begin{eqnarray} \label{1.4}
	2 \int_{\Omega}|\phi|^2 \log \left(\frac{|\phi|}{\|\phi\|_2}\right) d x+n(1+\log \mu)\|\phi\|_2^2 \leq \frac{a^2}{\pi}\|\nabla \phi\|_2^2, \quad (\text{see } \cite{DelPino2004,2001LL}).
\end{eqnarray}

In this paper, we aim to extend the main results of \cite{2025DZ} (for the case \(p=r=2\) of IBVP \eqref{1.3}) to the IBVP \eqref{1.1}. To treat the fractional equation, we use the fractional logarithmic Sobolev inequality (see \cref{l3.2}) established in \cite{2024CR}.
The explicit constant in the fractional logarithmic Sobolev inequality is 
$
a^2 K(n,s).
$
This constant is a function of \(n\) and \(s\), and it is much more complicated than the Lieb-Loss sharp constant \(\frac{a^2}{\pi}\) in \eqref{1.4}. This requires more delicate analysis and refined inequality arguments throughout our derivations and proofs. Recent works \cite{2021DWZ,2022DZ,2025DZ,2026FZ,2024LF,2024SHO} have developed new techniques for infinite time blow-up with high initial energy in strongly damped models. To our knowledge, the present work is among the first to prove such a result for the fractional Laplacian with strong damping. We provide explicit verifiable sufficient conditions and quantitatively reveal how fractional order (s) and dimension (n) shape the dynamical behavior of solutions.

The remainder of this paper is organized as follows. Section 2 states the main results. Section 3 presents several key lemmas. Section 4 provides the detailed proofs of the main results.

\section{Main results}

Throughout this paper, the following standard notations are adopted.
	$(\cdot,\cdot)$ represents the $L^2(\Omega)$-inner product and $\|\cdot\|_\gamma$ represents the $L^\gamma(\Omega)$-norm for $1 \leq \gamma \leq \infty$.
	$[\cdot]_{s}$ denotes the Gagliardo semi-norm, i.e.,
	$$
	[u]_s:=\left(\iint_Q \frac{|u(x)-u(y)|^2}{|x-y|^{n+2 s}} dx dy\right)^{\frac{1}{2}},
	$$
	where $Q=\mathbb{R}^{2 n} \setminus(\mathcal{C}(\Omega) \times \mathcal{C}(\Omega))$ and $\mathcal{C}(\Omega)=\mathbb{R}^n \setminus \Omega$. $W$ is a linear space of Lebesgue measurable functions from $\mathbb{R}^n$ to $\mathbb{R}$ such that the restriction to $\Omega$ of any function $u$ in $W$ belongs to $L^2(\Omega)$ and
	$$
	\iint_Q \frac{|u(x)-u(y)|^2}{|x-y|^{n+2 s}} dx dy<\infty.
	$$
	The space $W$ is equipped with the norm
	$$
	\|u\|_W=\left(\|u\|_2^2+\iint_Q \frac{|u(x)-u(y)|^2}{|x-y|^{n+2 s}} dx dy\right)^{\frac{1}{2}}.
	$$
    We shall work in the closed linear subspace
	$$
	W_0=\left\{u \in W: u(x)=0 \text { a.e. in } \mathbb{R}^n \backslash \Omega\right\}.
	$$
	An equivalent norm on $W_0$ is
	$$
	\|u\|_{W_0}^2:=[u]_s^2=\iint_Q \frac{|u(x)-u(y)|^2}{|x-y|^{n+2 s}} dx dy, \quad u \in W_0 .
	$$
	and the corresponding inner product on $W_0$ is defined by
	$$
	(u, v)_{W_0} = \iint_Q \frac{(u(x)-u(y))(v(x)-v(y))}{|x-y|^{n+2 s}} dx dy .
	$$
	Then
	$$
	\big((-\Delta)^{\frac{s}{2}} u,\ (-\Delta)^{\frac{s}{2}} v\big) = \frac{1}{2} C(n,\, s)\, (u,\, v)_{W_0}, \quad \forall\, u,\, v \in W_0.
	$$
	$$
    \left\|(-\Delta)^{\frac{s}{2}} u\right\|_{2}^2=\frac{1}{2} C(n, s)\|u\|_{W_0}^2, \quad \forall u \in W_0.
	$$
	$W^{-s}$ denotes the dual space of $W_0$, and $L^{\frac{2+\delta}{1+\delta}}(\Omega)$ denotes the dual space of $L^{2+\delta}(\Omega)$, where $\delta$ satisfies
	$
	0 < \delta < \frac{2s}{n-2s}.
	$
	$\langle \cdot,\cdot \rangle$ represents the duality bracket between $W_0$ and $W^{-s}$, and $\langle \cdot,\cdot \rangle_{**}$ represents the duality bracket between $L^{2+\delta}(\Omega)$ and $L^{\frac{2+\delta}{1+\delta}}(\Omega)$.
Clearly, the continuous embedding $W_0 \hookrightarrow L^{2+\delta}(\Omega)$ holds. In what follows, let $C_* = C_*(n,s)$ be the best constant for the embedding $W_0 \hookrightarrow L^{2+\delta}(\Omega)$. Let
$$\lambda_{1,s}= \lambda_{1,s}(n,s) =  \inf_{u\in W_0,\,u\neq 0} \frac{\left\|(-\Delta)^{\frac{s}{2}} u\right\|_{2}^2}{\|u\|_2^2} = \frac{1}{2} C(n, s) \inf_{u\in W_0,\,u\neq 0} \frac{\|u\|_{W_0}^2}{\|u\|_2^2}.$$ 
which is the first eigenvalue of the fractional Laplacian $(-\Delta)^s$ subject to homogeneous exterior Dirichlet boundary conditions. For further details concerning fractional Sobolev spaces, we refer readers to \cite{2012DNP}.
\begin{definition} \label{d2.1}
	Assume $0 < T \leq \infty$. Let $u_1 \in L^2(\Omega)$ and $u_0 \in W_0$. The function $u \in C([0,T);W_0)$ with
	$
	u_t \in L^2\left(0,T;W_0\right) \cap C\left([0,T);L^2(\Omega)\right), u_{tt} \in L^2\left(0,T;W^{-s}\right)
	$
	is called a weak solution to IBVP \eqref{1.1} over $[0,T)$, if $u(0)=u_0$, $u_t(0)=u_1$ and
\begin{eqnarray} \label{2.1}
		\left< u_{tt},\psi \right> + \frac{1}{2} C(n, s) (u_t,\psi)_{W_0} + \frac{1}{2} C(n, s) (u,\psi)_{W_0} = \left< u \ln|u|,\psi \right>_{**}
\end{eqnarray}
	holds for a.e. $t\in (0,T)$ and any $\psi\in W_0$.
\end{definition}

For any \( u\in W_0\), we define the following functionals
	\begin{align}
	J(u): & =\frac{1}{2} \left\|(-\Delta)^{\frac{s}{2}} u\right\|_{2}^2 + \frac{1}{4} \|u\|_2^2 - \frac{1}{2} \int_\Omega |u|^2 \ln|u| dx \nonumber\\
	&=\frac{1}{4} C(n, s) \|u\|_{W_0}^2 + \frac{1}{4} \|u\|_2^2 - \frac{1}{2} \int_\Omega |u|^2 \ln|u| dx, \label{2.2} \\
	I(u):&= \left\|(-\Delta)^{\frac{s}{2}} u\right\|_{2}^2  - \int_\Omega |u|^2 \ln|u| dx \nonumber \\
	&=\frac12 C(n, s) \|u\|_{W_0}^2 - \int_\Omega |u|^2 \ln|u| dx. \label{2.3}
	\end{align}
Then
\begin{eqnarray}  \label{2.4}
	J(u) = \frac{1}{4} \|u\|_2^2 + \frac{1}{2} I(u).
\end{eqnarray}
Let
$
d:=\inf_{u\in \mathcal{N}} J(u)
$
denote the mountain-pass level, where the Nehari manifold $\mathcal{N}$ is given by
$$
\mathcal{N} = \left\{ u\in W_0\setminus \{0\} \mid I(u)=0 \right\}.
$$
In addition, we introduce two auxiliary sets
\[
\mathcal{W}:=\left\{ u\in W_0 \mid I(u)>0 \right\} \cup \{0\},
\]
\[
\mathcal{V}:=\left\{ u\in W_0 \mid I(u)<0 \right\}.
\]
We also define the total energy functional as
\begin{eqnarray} \label{2.5}
	E(t) :=\frac{1}{2} \|u_t\|_2^2 + J(u).
\end{eqnarray}

First, the local existence theorem is presented. 

\begin{theorem} \label{t2.1}
	Assume $\widehat{\Lambda}>0$ is a constant. If $u_1\in L^2(\Omega)$ and $u_0\in W_0$ satisfy $\|u_0\|_{W_0}^2 + \|u_1\|_2^2 \leq \widehat{\Lambda}$, then there exists a constant $T_{\widehat{\Lambda}}>0$ only depending on $\widehat{\Lambda}$ such that IBVP \eqref{1.1} admits a weak solution $u$ defined on $[0,T_{\widehat{\Lambda}}]$, and the time interval $[0,T_{\widehat{\Lambda}}]$ can be maximally extended to $[0,T)$. In addition, the energy dissipation identity
	\begin{eqnarray} \label{2.6}
		E(t) + \frac{1}{2} C(n, s) \int_0^t \|u_\tau\|_{W_0}^2 d\tau = E(0)
	\end{eqnarray} 
	holds for all $t\in [0,T)$. Moreover, if $T<\infty$, then
	\begin{eqnarray} \label{2.7}
		\lim_{t\to T^-} \left( \|u_t(t)\|_2^2 + \frac{1}{2} C(n, s) \|u(t)\|_{W_0}^2 \right) = \infty.
	\end{eqnarray}
\end{theorem}
 
 Next, we have the following result on global existence.
 
\begin{theorem} \label{t2.2}
	Let $u_1\in L^2(\Omega)$ and $u_0\in W_0$. If $I(u_0)>0$ and $E(0)\leq d$, then IBVP \eqref{1.1} admits a global weak solution \(u\).
\end{theorem}

Finally, two results on infinite time blow-up are presented.

\begin{theorem} \label{t2.3}
	Assume $K(n,s) < 2e^{2+\frac{2s}{n}} $, where $K(n,s)$ is defined in \cref{l3.2}. Let $u_1\in L^2(\Omega)$ and $u_0\in W_0$.
	If \(E(0) \leq d\) and \(I(u_0) < 0\), then IBVP \eqref{1.1} has a weak solution \(u\) that blows up at infinite time, satisfying
	$
	\lim_{t\to \infty} \left\|(-\Delta)^{\frac{s}{2}}u(t)\right\|_2^2 = \infty.
	$
\end{theorem}

\begin{theorem} \label{t2.4}
	Assume $K(n,s) < e^2$, where $K(n,s)$ is defined in \cref{l3.2}. Let $u_1\in L^2(\Omega)$ and $u_0\in W_0$. Then
	
	(i) If
	\begin{eqnarray} \label{2.8}
		\begin{cases}
			I(u_0) < \|u_1\|_2^2,\\
			d < E(0) < \dfrac{\sqrt{\|(-\Delta)^{\frac{s}{2}}u_0\|_2^2 + 2(u_0,u_1)}}{2\sqrt{e b_*}},
		\end{cases}
	\end{eqnarray}
	where
	\begin{eqnarray} \label{2.9}
		b_* = \max\left\{ 2,\,\frac{1}{F_1(a_*)} \right\},
	\end{eqnarray}
	where $a_*>0$ is a constant satisfying (as shown in \Cref{r2.1})
	\[
	F_1(a_*) = \sup_{a\in\{a>0:F_2(a)\geq 0\}} F_1(a) \geq  \frac{\lambda_{1,s}\,s\Big(2-e^{-2} K(n,s)\Big) + s}{4s\big(\lambda_{1,s}+1\big)} > 0,
	\]
	and
	\begin{eqnarray} \label{2.10}
		F_1(a) = \frac{\lambda_{1,s}\,s\Big(2-a^2 K(n,s)\Big) + s + n\big(1+\ln a\big)}{4s\big(\lambda_{1,s}+1\big)},\quad
		F_2(a) = \frac{s\Big(1-a^2 K(n,s)\Big) - n\big(1+\ln a\big)}{4s\big(\lambda_{1,s}+1\big)},
	\end{eqnarray}
	then IBVP \eqref{1.1} has a weak solution $u$ that blows up at infinite time, satisfying $\lim_{t\to\infty} \|(-\Delta)^{\frac{s}{2}}u(t)\|_2^2 = \infty$.
	
	(ii) For any $\hat{Q} > d$, there exists a weak solution $u$ to IBVP \eqref{1.1} with initial data $u_0^{\hat{Q}} \in W_0$ and $u_1^{\hat{Q}} \in L^2(\Omega)$ satisfying the conditions in (i) and $E(0;(u_0^{\hat{Q}},u_1^{\hat{Q}})) = \hat{Q}$, which blows up at infinite time.
\end{theorem}

\begin{remark} \label{r2.1}
	The positivity of \(F_1(a_*)\) is verified here. Indeed, a straightforward calculation shows that \(F_1(a)\) is increasing on \((0, a_1)\) and decreasing on \((a_1, \infty)\), where $a_1 := \sqrt{\frac{n}{2\lambda_{1,s}\,s\,K(n,s)}}$. Recall $K(n,s) < e^2$, then 
	$$
	F_1(a_1) = \sup_{a\in(0,\infty)} F_1(a) \geq F_1\left(e^{-1}\right) = \frac{\lambda_{1,s}\,s\Big(2-e^{-2} K(n,s)\Big) + s}{4s\big(\lambda_{1,s}+1\big)} > 0.
	$$
	Since $\lim_{a \to 0^+} F_1(a) = \lim_{a \to +\infty} F_1(a) = -\infty$, there exist constants $0 < a_2 < a_1 < a_3$ such that $F_1(a) > 0$ if and only if $a \in (a_2, a_3)$. In addition, since $F_2(a)$ is strictly decreasing for all $a > 0$, and $\lim_{a \to 0^+} F_2(a) = \infty$, $ \lim_{a \to +\infty} F_2(a) = -\infty$, there exists a unique $a_4 > 0$ such that $F_2(a) > 0$ if and only if $a \in (0, a_4)$. Recall $K(n,s) < e^2$, then 
	$$
	F_2(e^{-1}) = \frac{s\Big(1-e^{-2} K(n,s)\Big) }{4s\big(\lambda_{1,s}+1\big)} > 0.
	$$
	Therefore, we have $F_1(e^{-1}) > 0$ and $F_2(e^{-1}) > 0$, which yields $e^{-1} \in I := (a_2, a_3) \cap (0, a_4)$, hence $I$ is nonempty. Moreover, it holds
	$$
	\sup_{a \in \{a>0:F_2(a)\geq0\}} F_1(a) = \sup_{a \in I} F_1(a) = F_1(a_*) \geq F_1(e^{-1}) = \frac{\lambda_{1,s}\,s\Big(2-e^{-2} K(n,s)\Big) + s}{4s\big(\lambda_{1,s}+1\big)} > 0,
	$$
	where $a_* = \min\{a_1, a_4\}$.
\end{remark}

\begin{remark}
	In the limiting sense as $s\to 1^-$ and $n\to\infty$, \cref{t2.4} are consistent with \cite[Theorem~4]{2025DZ}. Indeed, as $s\to 1^-$ and $n\to\infty$, we have $a^2 K(n,s)\to \frac{a^2}{\pi}$ by \Cref{r3.1} and $\left\|(-\Delta)^{\frac{s}{2}}u_0\right\|_2^2 \longrightarrow \|\nabla u_0\|_2^2$. Furthermore, we note that $\lim_{s\to 1^-}\lambda_{1,s}=\lambda_1$, where $\lambda_1$ stands for the first eigenvalue of the Laplace operator $-\Delta$ subject to homogeneous Dirichlet boundary conditions. Consequently, the coefficients in \eqref{2.10} converge to
	$$
	F_1(a) \longrightarrow \frac{\lambda_1\left(2 - \dfrac{a^2}{\pi}\right) + 1 + n\big(1+\ln a\big)}{4\big(\lambda_1+1\big)},
	\qquad
	F_2(a) \longrightarrow \frac{\left(1 - \dfrac{a^2}{\pi}\right) - n\big(1+\ln a\big)}{4\big(\lambda_1+1\big)},
	$$
	The condition $K(n,s) < e^2$ is automatically satisfied in the limit since $1/\pi < e^2$. 
\end{remark}

\begin{remark}
	From \cref{l3.7}, $I(u(t)) < \|u_t(t)\|_2^2 - \mathscr{O}$ holds for all $t \in [0,T)$ provided $(u_0,u_1) \in \mathcal{U}_{\mathscr{O}}$.
	Therefore, we impose $I(u_0) < \|u_1\|_2^2$ in \eqref{2.8}.
	The inequality
	$
	d < E(0) < \frac{\sqrt{\|(-\Delta)^{\frac{s}{2}}u_0\|_2^2 + 2(u_0,u_1)}}{2\sqrt{e b_*}}
	$
	given in \eqref{2.8} ensures invariance of the set $\mathcal{U}_{\mathscr{O}}$. The choice of $a_*$ in \Cref{r2.1} minimizes $b_*$.
	This enlarges the set of initial data satisfying \eqref{2.8}.
\end{remark}

\begin{figure}[htbp]
 	\centering
 	\includegraphics[width=0.7\textwidth]{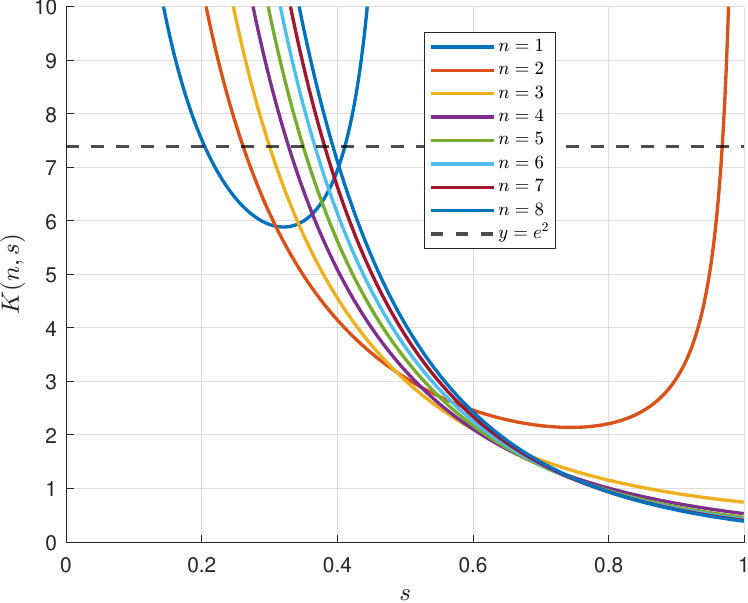}
 	\caption{Function $K(n,s)$ for $n=1,2,\dots,8$}
 	\label{fig:Kns}
 \end{figure}

\begin{remark}
	From \Cref{fig:Kns}, we can observe that the condition $K(n,s) < 2e^{2+\frac{2s}{n}}$ in Theorem 2.3 and the condition $K(n,s) < e^{2}$ in Theorem 2.4 are meaningful.
\end{remark}
 
\section{Preliminaries}

\begin{lemma} \label{l3.1}
	Let $u \in W_0 \setminus \{0\}$. Then
	
	(i) \( \lim_{\lambda\to 0} J(\lambda u) = 0 \) and \( \lim_{\lambda\to +\infty} J(\lambda u) = -\infty \);
	
	(ii) There exists a unique \( \lambda^* = \lambda^*(u) \) such that \( \frac{d}{d\lambda}J(\lambda u)|_{\lambda=\lambda^*} = 0 \).
	
	(iii) \( J(\lambda u) \) is increasing on \( (0,\lambda^*) \), decreasing on \( (\lambda^*,+\infty) \) and attains its maximum at \( \lambda^* \).
\end{lemma}

\begin{proof}
	(i) By \eqref{2.3} and \eqref{2.4}, it follows that
	$$
	\begin{aligned}
		J(\lambda u)&=\frac{\lambda^2}{4} C(n, s)\|u\|_{W_0}^2+\frac{\lambda^2}{4}\|u\|_2^2-\frac{1}{2}\int_\Omega|\lambda u|^2\ln|\lambda u|\,dx\\
		&=\frac{\lambda^2}{4}C(n, s)\|u\|_{W_0}^2+\frac{\lambda^2}{4}\|u\|_2^2-\frac{\lambda^2\ln\lambda}{2}\|u\|_2^2-\frac{\lambda^2}{2}\int_\Omega|u|^2\ln|u|\,dx \\
		& =\lambda^2\left(\frac{1}{4}C(n, s)\|u\|_{W_0}^2+\frac{1}{4}\|u\|_2^2-\frac{\ln\lambda}{2}\|u\|_2^2-\frac{1}{2}\int_\Omega|u|^2\ln|u|\,dx\right).
	\end{aligned}
	$$
	Then $\displaystyle\lim_{\lambda\to 0^+} J(\lambda u) = 0$ and $\displaystyle\lim_{\lambda\to +\infty} J(\lambda u) = -\infty$.
	
	(ii) Differentiate $J(\lambda u)$ with respect to $\lambda$ and simplify
	$$
	\begin{aligned}
		\frac{d}{d\lambda}J(\lambda u)&= \frac{\lambda}{2} C(n,s)\|u\|_{W_0}^2-\lambda\ln\lambda\|u\|_2^2-\lambda\int_\Omega|u|^2\ln|u|\,dx\\
		&=\frac{1}{\lambda}\left(\frac{\lambda^2}{2}C(n,s)\|u\|_{W_0}^2-\lambda^2\int_\Omega|u|^2\big(\ln\lambda+\ln|u|\big)\,dx\right)\\
		&=\frac{1}{\lambda}I(\lambda u).
	\end{aligned}
	$$
	Thus $\dfrac{d}{d\lambda}J(\lambda u)=0$ is equivalent to $I(\lambda u)=0$. Define an auxiliary function
	$$
	h(\lambda)=\frac{I(\lambda u)}{\lambda^2}
	=\frac12 C(n,s)\|u\|_{W_0}^2-\|u\|_2^2\ln\lambda-\int_\Omega|u|^2\ln|u|\,dx,
	$$
	so that $I(\lambda u)=0$ if and only if $h(\lambda)=0$. Since $u\in W_0\setminus\{0\}$, we have $\|u\|_2>0$. Note that
	$
	h'(\lambda)=-\frac{1}{\lambda}\|u\|_2^2<0, \forall \lambda>0.
	$
	Hence $h(\lambda)$ is strictly decreasing on $(0,+\infty)$. Combining the limits $\displaystyle\lim_{\lambda\to0^+}h(\lambda)=+\infty$ and $\displaystyle\lim_{\lambda\to+\infty}h(\lambda)=-\infty$, the intermediate value theorem guarantees a unique $\lambda^*>0$ satisfying $h(\lambda^*)=0$, i.e., $I(\lambda^* u)=0$. 
	
	(iii) Follows immediately from the sign of $\frac{d}{d\lambda} J(\lambda u)$ in (ii).
\end{proof}

\begin{lemma} \cite{2024CR} \label{l3.2}
	For any $u \in H^s(\mathbb{R}^n) \setminus \{0\}$ and $a > 0$, there holds
	$$
	2 \int_{\mathbb{R}^n} |u|^2 \ln \left( \frac{|u|}{\|u\|_2} \right) dx + \frac{n}{s} (1 + \ln a) \|u\|_2^2 \leq a^2 K(n,s) \left\| (-\Delta)^{s/2} u \right\|_2^2,
	$$
	where $K(n,s)=\frac{ne}{2s}\,C_{n,s}$ and
	$$C_{n,s}=\frac{\Gamma\left(\frac{n-2s}{2}\right)}{2^{2s}\pi^s\Gamma\left(\frac{n+2s}{2}\right)}
	\left(\frac{\Gamma(n)}{\Gamma(\frac n2)}\right)^{\frac{2s}{n}}.$$
	
	In particular, if $\Omega \subset \mathbb{R}^n$ is a bounded domain, then for any $u \in W_0 \setminus \{0\}$ and $a > 0$, there holds
	\begin{eqnarray} \label{3.1}
		2 \int_{\Omega} |u|^2 \ln \left( \frac{|u|}{\|u\|_2} \right) dx + \frac{n}{s} (1 + \ln a) \|u\|_2^2 \leq a^2 K(n,s) \left\| (-\Delta)^{s/2} u \right\|_2^2.
	\end{eqnarray}
\end{lemma}

\begin{remark} \label{r3.1}
	The explicit constant \( a^2 K(n,s) \) in the above inequality is not sharp for finite \(n\). In the special case \(s=1\) and $n\to\infty$, the asymptotic form of the estimate coincides with the sharp Lieb-Loss constant for the integer-order logarithmic Sobolev inequality.
\end{remark}


\begin{lemma} \label{l3.3}
	Let $u \in W_0 \setminus \{0\}$. If $I( u) \leq 0$, then $\|u\|_{W_0}\geq r_*:=\sup_{\varepsilon\in\left(0,\frac{2s}{n-2s}\right)}r_*(\varepsilon) $, where
	$ r_*(\varepsilon)=\left(\frac{\varepsilon\, C(n,s)}{2\, C_*^{2+\varepsilon}}\right)^{\frac{1}{\varepsilon}}, $
	and $C_* = C_*(n,s)$ denotes the best embedding constant for $W_0\hookrightarrow L^{2+\varepsilon}(\Omega)$.
\end{lemma}

\begin{proof}
	By $I(u)\leq0$ and the inequality $\ln|t|\leq\dfrac{|t|^\varepsilon}{\varepsilon}$ for $|t|\geq1$ and $\varepsilon>0$, we have
	$$
	\begin{aligned}
		\frac12 C(n,s)\|u\|_{W_0}^2 & \leq\int_\Omega|u|^2\ln|u| dx \\
		&=\int_{\{|u|\geq1\}}|u|^2\ln|u|\,dx+\int_{\{0<|u|<1\}}|u|^2\ln|u| dx\\
		&\leq\frac{1}{\varepsilon}\int_\Omega|u|^{2+\varepsilon} dx=\frac{1}{\varepsilon}\|u\|_{2+\varepsilon}^{2+\varepsilon} \leq \frac{1}{\varepsilon}C_*^{2+\varepsilon}\|u\|_{W_0}^{2+\varepsilon}.
	\end{aligned}
	$$
	Then
	$$
	\|u\|_{W_0}\geq \left(\frac{\varepsilon\, C(n,s)}{2\, C_*^{2+\varepsilon}}\right)^{\frac{1}{\varepsilon}}=r_*(\varepsilon).
	$$
	Hence $\|u\|_{W_0}\geq r_*>0$.
\end{proof}

\begin{lemma} \label{l3.4}
	The potential well depth
	$
	d \geq \frac{1}{4} e^{\frac{n}{s}(1+\ln a_0)} > 0,
	$
	where $a_0 = \sqrt{\dfrac{2}{K(n,s)}}>0$.
\end{lemma}

\begin{proof}
	Since $u \in \mathcal{N}$, then
	$$
	J(u) = \frac{1}{4}\|u\|_2^2.
	$$
	By \eqref{3.1}, for any $u \in W_0\setminus\{0\}$ and any $a>0$, we obtain
	$$
	2\int_\Omega |u|^2 \ln|u|\,dx \leq a^2 K(n,s) \left\|(-\Delta)^{\frac{s}{2}}u\right\|_2^2 + 2\|u\|_2^2 \ln\|u\|_2 - \frac{n}{s}(1+\ln a)\|u\|_2^2.
	$$
	Note that $\int_\Omega |u|^2 \ln|u|\,dx = \left\|(-\Delta)^{\frac{s}{2}} u\right\|_{2}^2$ from $I(u) = 0$,
	$$
	2\left\|(-\Delta)^{\frac{s}{2}} u\right\|_{2}^2 \leq a^2 K(n,s) \left\|(-\Delta)^{\frac{s}{2}}u\right\|_2^2 + 2\|u\|_2^2 \ln\|u\|_2 - \frac{n}{s}(1+\ln a)\|u\|_2^2.
	$$
	Then
	\begin{eqnarray} \label{3.2}
		\bigl(2 - a^2 K(n,s)\bigr) \left\|(-\Delta)^{\frac{s}{2}}u\right\|_2^2 \leq 2\|u\|_2^2 \ln\|u\|_2 - \frac{n}{s}(1+\ln a)\|u\|_2^2.
	\end{eqnarray}
	By \cref{l3.2}, $a^2 K(n,s)$ is strictly increasing and continuous with respect to $a>0$, and satisfies
	$$
	\lim_{a\to 0^+} a^2 K(n,s) = 0, \quad \lim_{a\to +\infty} a^2 K(n,s) = +\infty.
	$$
	By the intermediate value theorem, there exists a unique constant $a_0>0$ depending only on $n$ and $s$ such that $a_0^2 K(n,s) = 2$.	Substitute $a=a_0$ into \eqref{3.2},
	$$
	0 \leq 2\|u\|_2^2 \ln\|u\|_2 - \frac{n}{s}(1+\ln a_0)\|u\|_2^2.
	$$
	Since $u \in \mathcal{N}$, $\int_\Omega |u|^2 \ln|u|\,dx = \left\|(-\Delta)^{\frac{s}{2}} u\right\|_{2}^2 > 0$.
	Hence $\|u\|_2 > 0$ and 
	$$
	\|u\|_2^2 \geq e^{\frac{n}{s}(1+\ln a_0)}  > 0.
	$$
	which implies
	$
	d \geq \frac{1}{4} e^{\frac{n}{s}(1+\ln a_0)} > 0.
	$
\end{proof}

\begin{lemma} \label{l3.5}
	Let \(u\) be the weak solution of IBVP \eqref{1.1} with maximal existence time \([0,T)\). Then
	
	(i) If \(E(0) \leq  d\) and \(I(u_0)>0\), then \(I(u(t))>0\) for all \(t\in[0,T)\) (i.e., \(u(t)\in\mathcal{W}\)).
	
	(ii) If \(E(0) \leq d\) and \(I(u_0)<0\), then \(I(u(t))<0\) for all \(t\in[0,T)\) (i.e., \(u(t)\in\mathcal{V}\)).
\end{lemma}

\begin{proof}
	We only prove (ii), as the proof of (i) is similar. Arguing by contradiction, suppose there exists \(t_0 \in (0, T)\) such that \(I(u(t)) < 0\) for \(t \in [0, t_0)\) and \(I(u(t_0)) = 0\). Furthermore, \cref{l3.3} and the continuity of \(I(u)\) imply that \(u(t_0) \neq 0\). By the energy dissipation identity, we have
	\begin{eqnarray} \label{3.3}
		J(u(t)) \leq E(t) = E(0) - \frac{C(n,s)}{2}\int_{0}^{t} \|u_\tau\|_{W_0}^2\,d\tau \leq E(0), \quad \forall t \in [0, T).
	\end{eqnarray}
	
	If \(E(0) < d\), since \(u(t_0) \in \mathcal{N}\), we have \(J(u(t_0)) \geq d\) by definition of \(d\), which contradicts \(J(u(t_0)) \leq E(0) < d\) from \eqref{3.3}.
	
	If \(E(0) = d\), from \eqref{3.3}, we get \(J(u(t_0)) \leq d\). Together with \(J(u(t_0)) \geq d\), which holds since \(u(t_0) \in \mathcal{N}\), we conclude that \(J(u(t_0)) = d\). This equality implies
	$
	\dfrac{C(n,s)}{2}\displaystyle\int_{0}^{t_0} \|u_\tau\|_{W_0}^2\,d\tau = 0.
	$
	Hence \(\|u_\tau\|_{W_0} = 0\) for a.e. \(\tau\in[0,t_0]\). Since \(u_t(x) = 0\) for \(x\in \mathbb{R}^n\setminus \Omega\) and \(u_t \in C\left([0,T);L^2(\Omega)\right)\), it follows that \(u_t \equiv 0\) on \([0,t_0]\), so that \(u(t) \equiv u_0\) for all \(t\in[0,t_0]\). Consequently, \(I(u_0) = I(u(t_0)) = 0\), which contradicts the assumption \(I(u_0) < 0\). This completes the proof.
\end{proof}

\begin{lemma} \label{l3.6}
	For any \( u \in W_0 \setminus \{0\}\), if \(I( u) < 0\), it follows that
	\begin{eqnarray} \label{3.4}
		I( u) < 2\left(J( u) - d\right).
	\end{eqnarray}
\end{lemma}

\begin{proof}
	\cref{l3.1} implies there exists \(\alpha^* \in (0,1)\) such that \(I(\alpha^*  u) = 0\), i.e., \(\alpha^*  u \in \mathcal{N}\). By the definition of \(d\), we have \(J(\alpha^*  u) \geq d\).
	Since \(u\in W_0\setminus\{0\}\), we get \(\|u\|_2>0\).
	Define
	\[g(\alpha) := 2 J(\alpha  u) - I(\alpha  u) = \frac{\alpha^2}{2} \| u\|_2^2, \quad \alpha>0.\]
	Differentiating \(g(\alpha)\) with respect to \(\alpha\) yields
	\begin{align*}
		g'(\alpha) = \alpha \| u\|_2^2 > 0,\quad \forall\alpha>0.
	\end{align*}
	Thus \(g(\alpha)\) is strictly increasing on \((0,+\infty)\). Then
	\[ 2 J( u) - I( u) = g(1) > g(\alpha^*) = 2 J(\alpha^*  u) - I(\alpha^*  u) = 2 J(\alpha^*  u) \geq 2 d.\]
	This completes the proof.
\end{proof}

\begin{lemma} \label{l3.7}
	Assume $K(n,s) < 2e^2$. For any fixed $\mathscr{O} > 0$, define
	\begin{eqnarray} \label{3.5}
		\mathcal{U}_{\mathscr{O}} := \Big\{ (\Phi, \Psi) \in W_0 \times L^2(\Omega) :\  I(\Phi) < \|\Psi\|_2^2 - \mathscr{O},  d < E(\Phi, \Psi) \leq \frac{\sqrt{\|(-\Delta)^{\frac{s}{2}}\Phi\|_2^2 + 2(\Phi, \Psi)}}{2\sqrt{e b_*}} - \frac{\mathscr{O}}{2} \Big\},
	\end{eqnarray}
	where $b_* > 0$ is given by \eqref{2.9}, $I(\Phi)=\|(-\Delta)^{\frac{s}{2}}\Phi\|_2^2 - \int_\Omega |\Phi|^2 \ln|\Phi| dx$ and 
	$
	E(\Phi, \Psi) = \frac{1}{2}\|\Psi\|_2^2 + \frac{1}{2}\|(-\Delta)^{\frac{s}{2}}\Phi\|_2^2 + \frac{1}{4}\|\Phi\|_2^2 - \frac{1}{2}\int_\Omega \Phi^2 \ln|\Phi| \, dx.
	$
	If $(u_0, u_1) \in \mathcal{U}_{\mathscr{O}}$, then $I(u(t)) < \|u_t(t)\|_2^2 - \mathscr{O}$ for all $t \in [0, T)$, and
	\begin{eqnarray} \label{3.6}
		d < E(0) \leq \frac{\sqrt{\|(-\Delta)^{\frac{s}{2}}u(t)\|_2^2 + 2(u(t), u_t(t))}}{2\sqrt{e b_*}} - \frac{\mathscr{O}}{2}, \quad \forall t \in [0, T).
	\end{eqnarray}
\end{lemma}

\begin{proof}
	Define $\mathfrak{A}(u(t)) = I(u(t)) - \|u_t(t)\|_2^2$ for $t \in [0, T)$. 
	We first prove that $\mathfrak{A}(t) < -\mathscr{O}$ for all $t \in [0, T)$. Suppose otherwise. Since the initial condition $(u_0, u_1) \in \mathcal{U}_{\mathscr{O}}$ yields $\mathfrak{A}(0) = I(u_0) - \|u_1\|_2^2 < -\mathscr{O}$, there exists $t_0 \in (0, T)$ such that $\mathfrak{A}(t_0) = -\mathscr{O}$ and $\mathfrak{A}(t) < -\mathscr{O}$ for all $t \in [0, t_0)$. Let
	\[
	\mathcal{E}(t) = \|u(t)\|_2^2 + \int_0^t \|(-\Delta)^{\frac{s}{2}}u(\tau)\|_2^2 \, d\tau, \quad t \in [0, T).
	\]
	Direct computation using \eqref{1.1} and \eqref{2.3} yields
	\[
	\mathcal{E}'(t) = 2(u(t), u_t(t)) + \|(-\Delta)^{\frac{s}{2}}u(t)\|_2^2,
	\]
	\[
	\mathcal{E}''(t) = 2\|u_t\|_2^2 - 2I(u(t)) = -2\mathfrak{A}(u(t)) > 0, \quad t \in [0, t_0).
	\]
	Thus $\mathcal{E}'(t)$ is strictly increasing on $[0, t_0]$. Since $(u_0, u_1) \in \mathcal{U}_{\mathscr{O}}$, we obtain
	\begin{eqnarray} \label{3.7}
		\|(-\Delta)^{\frac{s}{2}}u(t_0)\|_2^2 + 2(u(t_0), u_t(t_0)) = \mathcal{E}'(t_0) \geq \mathcal{E}'(0) = \|(-\Delta)^{\frac{s}{2}}u_0\|_2^2 + 2(u_0, u_1) \geq 4 e b_* \left( \frac{\mathscr{O}}{2} + E(0) \right)^2. 
	\end{eqnarray}
	Moreover, by \eqref{2.6}, $\mathfrak{A}(t_0) = -\mathscr{O}$ together with \eqref{2.4} and \eqref{2.5} yields
	\[
	\begin{aligned}
		E(0) \geq E(t_0) &= \frac{1}{2}\|u_t(t_0)\|_2^2 + J(u(t_0)) \\
		&= \frac{1}{2}\|u_t(t_0)\|_2^2 + \frac{1}{2}I(u(t_0)) + \frac{1}{4}\|u(t_0)\|_2^2 \\
		&= \|u_t(t_0)\|_2^2 - \frac{\mathscr{O}}{2} + \frac{1}{4}\|u(t_0)\|_2^2 \\
		&\geq \frac{1}{4}\|u(t_0)\|_2^2 - \frac{\mathscr{O}}{2}.
	\end{aligned}
	\]
	This yields $4E(0) + 2\mathscr{O} \geq \|u(t_0)\|_2^2$. By \cref{l3.4}, $d \geq \frac{1}{4} e^{\frac{n}{s}(1+\ln a_0)}$ with $a_0 = \sqrt{\frac{2}{K(n,s)}}>0$. Since $K(n,s) < 2e^2$, it follows that $d > \frac14$. Hence $4E(0) + 2\mathscr{O} > 4d > 1$, so
	\begin{eqnarray} \label{3.8}
		\|u(t_0)\|_2^2 \ln\|u(t_0)\|_2^2 \leq (4E(0) + 2\mathscr{O}) \ln(4E(0) + 2\mathscr{O}).
	\end{eqnarray}
	By \eqref{3.1}, for any $u \in W_0\setminus\{0\}$ and any $a>0$, we obtain
	$$
	2\int_\Omega |u|^2 \ln|u|\,dx \leq a^2 K(n,s)\left\|(-\Delta)^{\frac{s}{2}}u\right\|_2^2 + 2\|u\|_2^2 \ln\|u\|_2 - \frac{n}{s}(1+\ln a)\|u\|_2^2,
	$$
	which, together with \eqref{2.2}, \eqref{2.5}, \eqref{2.6} and \eqref{3.8}, implies
	\begin{eqnarray} \label{3.9}
		\begin{aligned}
			E(0) \geq E(t_0)
			& = \frac12\|u_t(t_0)\|_2^2 + \frac12\left\|(-\Delta)^{\frac{s}{2}}u(t_0)\right\|_2^2 - \frac12\int_\Omega \big|u(t_0)\big|^2\ln\big|u(t_0)\big|\,dx + \frac14\big\|u(t_0)\big\|_2^2 \\
			&\geq \frac12\big\|u_t(t_0)\big\|_2^2
			+ A(a) \left\|(-\Delta)^{\frac{s}{2}}u(t_0)\right\|_2^2
			+ B(a) \big\|u(t_0)\big\|_2^2 - \left( \frac{\mathscr{O}}{2} + E(0) \right) \ln(2\mathscr{O} + 4E(0)).
		\end{aligned}
	\end{eqnarray}
	where $A(a) = \frac{2 - a^2 K(n,s)}{4},\ B(a) =\frac14 + \frac{n(1+\ln a)}{4s}$. It is not difficult to verify that
	\[
	A(a) = F_1(a) + F_2(a),\quad F_1(a) = B(a) + \lambda_{1,s} F_2(a),
	\]
	where $F_1(a)$ and $F_2(a)$ are given by \eqref{2.10}. Note that $\lambda_{1,s}\|u\|_2^2 \le \left\|(-\Delta)^{\frac{s}{2}}u\right\|_2^2$. For any $a > 0$ satisfying $F_2(a) \geq 0$, 
	\[
	\begin{aligned}
		E(0) &\geq \frac12\big\|u_t(t_0)\big\|_2^2
		+ F_1(a) \left\|(-\Delta)^{\frac{s}{2}}u(t_0)\right\|_2^2 + \lambda_{1,s} F_2(a) \big\|u(t_0)\big\|_2^2 +
		B(a) \big\|u(t_0)\big\|_2^2 - \left( \frac{\mathscr{O}}{2} + E(0) \right) \ln(2\mathscr{O} + 4E(0)) \\
		&\geq \frac12\big\|u_t(t_0)\big\|_2^2
		+ F_1(a) \left\|(-\Delta)^{\frac{s}{2}}u(t_0)\right\|_2^2 + F_1(a) \big\|u(t_0)\big\|_2^2 - \left( \frac{\mathscr{O}}{2} + E(0) \right) \ln(2\mathscr{O} + 4E(0)).
	\end{aligned}
	\]
	In view of \Cref{r2.1}, we set $a = a_*$ in the preceding inequality. Together with \eqref{2.9} and \eqref{3.7}, this yields
	\[
	\begin{aligned}
		\frac{\mathscr{O}}{2} + E(0)
		&\geq \frac{1}{b_*}\left(\left\|(-\Delta)^{\frac{s}{2}}u(t_0)\right\|_2^2 + 2(u(t_0), u_t(t_0))\right) - \left( \frac{\mathscr{O}}{2} + E(0) \right) \ln(2\mathscr{O} + 4E(0)) \\
		&\geq 4e \left( \frac{\mathscr{O}}{2} + E(0) \right)^2 - \left( \frac{\mathscr{O}}{2} + E(0) \right) \ln(2\mathscr{O} + 4E(0)).
	\end{aligned}
	\]
	This leads to $2e \mathscr{O} + 4e E(0) \leq \ln(2e \mathscr{O} + 4e E(0))$, a contradiction in view of $x > \ln x$ for all $x > 0$.
	Consequently, for every $t \in [0, T)$, we have $\mathfrak{A}(u(t)) < -\mathscr{O}$.
	
	Next, we establish \eqref{3.6}. From the above, $\mathfrak{A}(t) < -\mathscr{O}$ for all $t \in [0, T)$, which implies
	$
	\mathcal{E}''(t) = -2\mathfrak{A}(t) > 0.
	$
	This yields
	\begin{align*}
		\left\|(-\Delta)^{\frac{s}{2}}u(t)\right\|_2^2 + 2(u(t), u_t(t)) = \mathcal{E}'(t) \geq \mathcal{E}'(0) = \left\|(-\Delta)^{\frac{s}{2}}u_0\right\|_2^2 + 2(u_0, u_1) \geq 4 e b_* \left( \frac{\mathscr{O}}{2} + E(0) \right)^2, \quad t \in [0, T).
	\end{align*}
	This completes the proof of the invariant set.
\end{proof}

\section{Proof of the Theorems}

{\bfseries Proof of \cref{t2.1}}

Let $\{e_i\}_{i=1}^{\infty}$ be the eigenfunctions of the fractional Laplacian $(-\Delta)^s$ subject to the homogeneous exterior Dirichlet boundary condition. Then $\{e_i\}_{i=1}^{\infty}$ forms an orthogonal basis of $W_0$. Define the Galerkin approximate solutions of IBVP \eqref{1.1} as
\[
u_k(x,t) = \sum_{i=1}^k h_{ik}(t) e_i(x), \quad k = 1,2,3,\cdots,
\]
which satisfy, for each $i = 1, 2, \cdots, k$,
\begin{eqnarray}\label{4.1}
	\begin{aligned}
		\langle u_{ktt}, e_i \rangle + \frac{1}{2}C(n,s)(u_{kt},e_i)_{W_0} + \frac{1}{2}C(n,s)(u_k,e_i)_{W_0}
		= \langle u_k \ln|u_k|, e_i \rangle_{**},
	\end{aligned}
\end{eqnarray}
with initial conditions $u_k(0) = u_{k0}$, $u_{kt}(0) = u_{k1}$, where $u_{k0}, u_{k1} \in \operatorname{span}\{e_1,e_2,\cdots,e_k\}$ are chosen such that
\begin{align}
	&u_{k0} \to u_0 \quad \text{in } W_0, \quad \text{as } k\to\infty, \label{4.2} \\
	&u_{k1} \to u_1 \quad \text{in } L^2(\Omega), \quad \text{as } k\to\infty. \label{4.3}
\end{align}
By the standard theory of ordinary differential equations, system \eqref{4.1} admits a unique local solution $u_k(t)$ on $[0, T_k]$ for some $T_k > 0$. Multiplying \eqref{4.1} by $h_{ik}'(t)$, summing over $i = 1,2,\cdots,k$, and integrating with respect to $t$ from $0$ to $t$, we obtain
\begin{align*}
	&\frac{1}{2}\|u_{kt}(t)\|_2^2 + \frac{1}{4}C(n,s)\|u_k(t)\|_{W_0}^2 + \frac{1}{2}C(n,s)\int_0^t \|u_{k\tau}(\tau)\|_{W_0}^2 d\tau \\
	&= \frac{1}{2}\|u_{k1}\|_2^2 + \frac{1}{4}C(n,s)\|u_{k0}\|_{W_0}^2 + \int_0^t \int_\Omega u_k(\tau) \ln|u_k(\tau)| \cdot u_{k\tau}(\tau) dx d\tau.
\end{align*}
Following the same argument as in \cite[(4.5)--(4.13)]{2025DZ}, we can show that there exists a constant $T_3 > 0$ depending only on $\widehat{\Lambda}$ such that
\begin{eqnarray} \label{4.4}
	\|u_{kt}(t)\|_2^2 + C(n,s)\|u_k(t)\|_{W_0}^2 + C(n,s)\int_0^t \|u_{k\tau}(\tau)\|_{W_0}^2 d\tau \leq \Lambda_1, \quad \forall t\in[0, T_3],
\end{eqnarray}
for some constant $\Lambda_1 > 0$ independent of $k$.

Therefore, there exists a function $u$ and a subsequence of $\{u_k\}_{k=1}^\infty$ (still denoted by $\{u_k\}$) such that
\begin{align}
	&u_k \stackrel{w^*}{\rightharpoonup} u \quad \text{in } L^\infty(0, T_3; W_0), \quad \text{as } k\to\infty, \label{4.5} \\
	&u_k \stackrel{w}{\rightharpoonup} u \quad \text{in } L^2(0, T_3; W_0), \quad \text{as } k\to\infty, \label{4.6} \\
	&u_{kt} \stackrel{w^*}{\rightharpoonup} u_t \quad \text{in } L^\infty(0, T_3; L^2(\Omega)), \quad \text{as } k\to\infty, \label{4.7} \\
	&u_{kt} \stackrel{w}{\rightharpoonup} u_t \quad \text{in } L^2(0, T_3; W_0), \quad \text{as } k\to\infty. \label{4.8}
\end{align}

By Lions' compactness theorem (see \cite{Simon1987}), \eqref{4.5} and \eqref{4.7} imply that
\[
u_k \to u \quad \text{strongly in } L^2(0, T_3; L^2(\Omega)), \quad \text{as } k\to\infty,
\]
and hence $u_k \to u$ almost everywhere in $\Omega\times[0, T_3]$. Therefore,
\begin{eqnarray}\label{4.9}
	u_k \ln|u_k| \to u \ln|u| \quad \text{a.e. in } \Omega\times[0, T_3], \quad \text{as } k\to\infty.
\end{eqnarray}

Furthermore, similar to the estimate in \cite[(4.23)]{2025DZ}, we deduce from \eqref{4.4} that
\begin{eqnarray} \label{4.10}
	\int_\Omega \left| u_k \ln|u_k| \right|^{\frac{2+\delta}{1+\delta}} dx \leq \left( \frac{1}{e} \right)^{\frac{2+\delta}{1+\delta}} |\Omega| + \frac{C_*^{2+\delta}}{\delta^{\frac{2+\delta}{1+\delta}}} \big(2\Lambda_1\big)^{\frac{2+\delta}{2}} =: \tilde{C},
\end{eqnarray}
where $C_* = C_*(n,s) > 0$ is the best constant of the embedding $H_0^{s}(\Omega) \hookrightarrow L^{2+\delta}(\Omega)$. Combining \eqref{4.10} with \eqref{4.9}, we conclude that there exists a subsequence such that
\begin{align}
	&u_k \ln|u_k| \stackrel{w^*}{\rightharpoonup} u \ln|u| \quad \text{in } L^{\infty}\left(0, T_3; L^{\frac{2+\delta}{1+\delta}}(\Omega)\right), \quad \text{as } k\to\infty;  \notag\\
	&u_k \ln|u_k| \stackrel{w}{\rightharpoonup} u \ln|u| \quad \text{in } L^2\left(0, T_3; L^{\frac{2+\delta}{1+\delta}}(\Omega)\right), \quad \text{as } k\to\infty. \label{4.11}
\end{align}

Let $P = \left\{ \sum_{i=1}^k \zeta_i e_i : \zeta_i\in\mathbb{R},\ k\in\mathbb{N} \right\}$. Then $P$ is dense in $W_0$. For any $\psi\in W_0$, there exists a sequence $\{\psi_i\}_{i=1}^\infty \subset P$ such that $\lim_{i\to\infty} \psi_i = \psi$ in $W_0$. By the construction of Galerkin approximations, \eqref{4.1} holds with $e_i$ replaced by $\psi_i$. Taking the limit $i \to \infty$, we obtain
\begin{align*}
	\left| \langle u_{ktt}, \psi \rangle \right|
	&\leq \frac{1}{2}C(n,s)\left| \left( u_{kt}, \psi \right)_{W_0} \right| + \frac{1}{2}C(n,s)\left| \left( u_k, \psi \right)_{W_0} \right| + \left| \langle u_k \ln|u_k|, \psi \rangle_{**} \right| \\
	&\leq \frac{1}{2}C(n,s)\|u_{kt}\|_{W_0}\|\psi\|_{W_0} + \frac{1}{2}C(n,s)\|u_k\|_{W_0}\|\psi\|_{W_0} + \| u_k \ln|u_k| \|_{\frac{2+\delta}{1+\delta}} \|\psi\|_{L^{2+\delta}(\Omega)}.
\end{align*}
Since $W_0\hookrightarrow L^{2+\delta}(\Omega)$, there exists a constant $C_* = C_*(n,s) > 0$ such that $\|\psi\|_{L^{2+\delta}(\Omega)} \leq C_* \|\psi\|_{W_0}$ for all $\psi\in W_0$. For any $\psi\in W_0$,
\[
\left| \langle u_{ktt}, \psi \rangle \right| \leq \widetilde{C_1} \|\psi\|_{W_0} \left( \|u_{kt}\|_{W_0} + \|u_k\|_{W_0} + \| u_k \ln|u_k| \|_{\frac{2+\delta}{1+\delta}} \right).
\]
From \eqref{4.10}, we immediately obtain
$
\|u_{ktt}\|_{W^{-s}} \leq \widetilde{C_1} \left( \|u_{kt}\|_{W_0} + \|u_k\|_{W_0} + \tilde{C}^{\frac{1+\delta}{2+\delta}} \right).
$
From the uniform bound \eqref{4.4}, we know $\|u_k\|_{W_0} \le \sqrt{2\Lambda_1/C(n,s)}$ holds uniformly for all $k$. From \eqref{4.8}, we immediately obtain
\begin{eqnarray}\label{4.12}
	u_{ktt} \stackrel{w}{\rightharpoonup} u_{tt} \quad \text{in } L^2(0, T_3; W^{-s}), \quad \text{as } k\to\infty.
\end{eqnarray}
Integrating \eqref{4.1} with respect to $t$ from $0$ to $t$ and taking the limit $k \to \infty$ via \eqref{4.8}, \eqref{4.11} and \eqref{4.12}, we conclude that for all $t \in [0, T_3]$,
\begin{align*}
	&(u_t(t), e_i) + \frac{1}{2}C(n,s)\int_0^t \left( u(\tau), e_i \right)_{W_0} d\tau + \frac{1}{2}C(n,s)\int_0^t \left( u_\tau(\tau), e_i \right)_{W_0} d\tau \\
	&= (u_1, e_i) + \int_0^t \langle u(\tau) \ln|u(\tau)|, e_i \rangle_{**} d\tau,
\end{align*}
for all $i = 1, 2, \cdots$. Since $\{e_i\}$ is dense in $W_0$, the weak formulation in \Cref{d2.1} holds for all test functions $\psi \in W_0$. Moreover, from \eqref{4.5} and \eqref{4.8}, we have $u \in C([0, T_3]; W_0)$. From \eqref{4.8}, \eqref{4.12} and \cite[Proposition 1.2]{2013SS}, we have $u_t \in C([0,T_3];L^2(\Omega))$. Combined with \eqref{4.2} and \eqref{4.3}, this yields $u(0) = u_0$ and $u_t(0) = u_1$. The energy identity \eqref{2.6} follows from multiplying the weak formulation by $u_t$ and integrating in time. Finally, we prove the blow-up criterion \eqref{2.7}. Suppose, for contradiction, that $T < \infty$ but
\[
\limsup_{t\to T^-} \big( \|u_t(t)\|_2^2 + \|u(t)\|_{W_0}^2 \big) < \infty.
\]
Then there exists a constant $R_* > 0$ such that
\[
\|u_t(t)\|_2^2 + \|u(t)\|_{W_0}^2 \leq R_*, \quad \forall t\in[0, T).
\]
By the local existence result, there exists $T_{R_*} > 0$ depending only on $R_*$ such that for any $t_0 \in [0, T)$, the solution can be extended to $[t_0, t_0 + T_{R_*}]$. Taking $t_0 = T - T_{R_*}/2$, we see that the solution exists on $[0, T + T_{R_*}/2]$, which contradicts the maximality of $T$. This completes the proof.
\hfill$\square$

{\bfseries Proof of \cref{t2.2}}

Only the case of \(E(0) < d\) is proved here, and the result for \(E(0) = d\) can be obtained by the same argument as in \cite[Theorems~2.2]{2025DZ}. By \cref{l3.5}, we have \(u(t) \in \mathcal{W}\) and thus \(I(u(t)) \geq 0\) for all \(t \in [0, T)\). In accordance with \eqref{2.4}, \eqref{2.5} and \eqref{2.6}, one has
\begin{align*}
	d > E(0) &= \frac{1}{2}C(n,s)\int_0^t \|u_\tau\|_{W_0}^2 d\tau + \frac{1}{2}\|u_t\|_2^2 + \frac{1}{4} \|u\|_2^2 + \frac{1}{2} I(u) \\
	&\geq \frac{1}{2}C(n,s)\int_0^t \|u_\tau\|_{W_0}^2 d\tau + \frac{1}{2}\|u_t\|_2^2,\quad t \in [0,T).
\end{align*}
Then we arrive at $\displaystyle\int_0^t \|u_\tau\|_{W_0}^2 d\tau < \dfrac{2d}{C(n,s)}$ and $\|u_t\|_2^2 < 2d$, which yields
\begin{align*}
	\|u\|_{W_0} &= \left\| u_0 + \int_0^t u_\tau d\tau \right\|_{W_0} \leq \|u_0\|_{W_0}+\int_0^t \|u_\tau\|_{W_0}\,d\tau \leq \|u_0\|_{W_0} + t^{\frac{1}{2}} \left(\int_0^t \|u_\tau\|_{W_0}^2 d\tau\right)^{\frac{1}{2}} \\
	&\leq \|u_0\|_{W_0} + t^{\frac{1}{2}} \left(\dfrac{2d}{C(n,s)}\right)^{\frac{1}{2}},\quad t \in [0,T).
\end{align*}
Thus, it holds
$
\|u_t\|_2^2 + \|u\|_{W_0}^2 \leq 2d + \left(\|u_0\|_{W_0} + t^{\frac{1}{2}} \left(\dfrac{2d}{C(n,s)}\right)^{\frac{1}{2}}\right)^2,\quad t \in [0,T).
$
By \cref{t2.1}, we know $T = \infty$.
\hfill$\square$

{\bfseries Proof of \cref{t2.3}}

Step 1: Extension of \(u\) to infinite time.

By \cref{l3.5}, we have \(u(t) \in \mathcal{V}\) and thus \(I(u(t)) < 0\) for all \(t \in [0, T)\).
Let $H(t) = \| u(t) \|_2^2$. A direct computation gives
\begin{eqnarray} \label{4.13}
	H'(t) = 2(u(t), u_t(t)) \leq \| u(t) \|_2^2 + \| u_t(t) \|_2^2.
\end{eqnarray}
Write $M(n,s) = -\frac{n}{2s} - \frac{n}{4s}\ln\frac{2}{K(n,s)} + 1$. In accordance with \eqref{4.13}, \eqref{2.5}, \eqref{2.6} and \eqref{3.4}, we deduce
\begin{eqnarray} \label{4.14}
	\begin{aligned}
		&(H(t) + e^{2M(n,s)})\ln(H(t) + e^{2M(n,s)}) - H'(t) \\
		&\geq \| u \|_2^2 \ln(\| u \|_2^2 + e^{2M(n,s)}) - \| u_t \|_2^2 - \| u \|_2^2 \\
		&\geq \| u \|_2^2 \ln(\| u \|_2^2 + e^{2M(n,s)}) + 2(J(u) - d) - \| u \|_2^2 \\
		&> \| u \|_2^2 \ln(\| u \|_2^2 + e^{2M(n,s)}) + I(u) - \| u \|_2^2 \\
		&= \| u \|_2^2 \ln(\| u \|_2^2 + e^{2M(n,s)}) + \left\|(-\Delta)^{\frac{s}{2}} u\right\|_{2}^2 - \int_\Omega |u|^2 \ln|u| \, dx - \| u \|_2^2.
	\end{aligned}
\end{eqnarray}
By \eqref{3.1} with $a = \sqrt{2/K(n,s)}$, we get
\[
\int_\Omega |u|^2 \ln|u| \, dx \leq \left\|(-\Delta)^{\frac{s}{2}} u\right\|_{2}^2 + \frac12 \| u \|_2^2 \ln\| u \|_2^2 - \frac{n}{2s}\left(1 + \frac{1}{2}\ln\frac{2}{K(n,s)}\right) \| u \|_2^2.
\]
Substituting this into \eqref{4.14},
\[
\begin{aligned}
	&(H(t) + e^{2M(n,s)})\ln(H(t) + e^{2M(n,s)}) - H'(t) \\
	&\geq \| u \|_2^2 \left( \frac 12 \ln(\| u \|_2^2 + e^{2M(n,s)}) - 1 + \frac{n}{2s} + \frac{n}{4s}\ln\frac{2}{K(n,s)} \right) \\
	&\geq \| u \|_2^2 \left( M(n,s) - 1 + \frac{n}{2s} + \frac{n}{4s}\ln\frac{2}{K(n,s)} \right) \geq 0.
\end{aligned}
\]
Consequently, we obtain
$
(H(t) + e^{2M(n,s)})\ln(H(t) + e^{2M(n,s)}) \geq H'(t).
$
By Gronwall's inequality,
\begin{eqnarray} \label{4.15}
	H(t) \leq (H(0) + e^{2M(n,s)})^{e^t} - e^{2M(n,s)} = (\| u_0 \|_2^2 + e^{2M(n,s)})^{e^t} - e^{2M(n,s)} : = H_1(t).
\end{eqnarray}
By an argument similar to that for \eqref{3.9}, combined with \eqref{4.15}, we obtain
\[
\begin{aligned}
	|E(0)| \geq E(t)
	& = \frac12\|u_t(t)\|_2^2 + \frac12 \left\|(-\Delta)^{\frac{s}{2}}u(t)\right\|_2^2 - \frac12\int_\Omega \big|u(t)\big|^2\ln\big|u(t)\big|\,dx + \frac14\big\|u(t)\big\|_2^2 \\
	&=  \frac12\big\|u_t(t)\big\|_2^2
	+ A(a) \left\|(-\Delta)^{\frac{s}{2}}u(t)\right\|_2^2
	+ B(a) \big\|u(t)\big\|_2^2 - \frac14 \|u(t)\|_2^2 \ln\|u(t)\|_2^2 \\
	&\geq \frac12\big\|u_t(t)\big\|_2^2
	+ A(a) \left\|(-\Delta)^{\frac{s}{2}}u(t)\right\|_2^2
	+ B(a) \big\|u(t)\big\|_2^2 - \frac14 (H_1(t)+1)\ln(H_1(t)+1).
\end{aligned}
\]
where $A(a) = \frac{2 - a^2 K(n,s)}{4}, B(a) =\frac14 + \frac{n(1+\ln a)}{4s}$. To guarantee the existence of some \(a>0\) such that both \(A(a)>0\) and \(B(a) \geq 0\), the interval 
$
\left[ e^{-1-\frac{s}{n}}
,\sqrt{\frac{2}{K(n,s)}}
\right)
$ must be non-empty, which implies that $K(n,s) <  2e^{2+\frac{2s}{n}}$. Hence,
\[
\left\|(-\Delta)^{\frac{s}{2}} u\right\|_{2}^2 + \| u_t \|_2^2 \leq \max \left\{2, \frac{1}{A(a)}\right\} \left(|E(0)| + \frac14 (H_1(t)+1)\ln(H_1(t)+1)\right).
\]
Let $\varpi(t) = \max \left\{2, \frac{1}{A(a)}\right\} \left(|E(0)| + \frac14 (H_1(t)+1)\ln(H_1(t)+1)\right)$, which is continuous on $[0, \infty)$. Since $\|u(t)\|_{W_0}^2 = \dfrac{2}{C(n,s)}\left\|(-\Delta)^{\frac{s}{2}}u(t)\right\|_2^2$, the above bound implies $\limsup_{t\to T^-}(\|u_t(t)\|_2^2 + \|u(t)\|_{W_0}^2) < \infty$. By the blow-up criterion in \cref{t2.1}, the maximal existence time satisfies $T = \infty$.

Step 2: Blow-up at infinite time. Let
\[
\mathscr{F}(t) := \| u(t) \|_2^2 + \frac{1}{2}C(n,s)\int_0^t \|u(\tau)\|_{W_0}^2 \, d\tau, \quad t \geq 0.
\]
By \eqref{1.1}, \eqref{2.3} and \(I(u) < 0\) for all \(t \in [0,T)\), for $t \geq 0$,
\begin{eqnarray} \label{4.16}
	\mathscr{F}'(t) = 2(u(t), u_t(t)) + \frac{1}{2}C(n,s)\|u(t)\|_{W_0}^2, \quad \mathscr{F}''(t) = 2\|u_t\|_2^2 -2I(u) > 0.
\end{eqnarray}

Case (I): \(E(0) < d\). It follows from \eqref{2.6}, \eqref{4.16} and \eqref{3.4} that for $t \geq 0$,
\[
\begin{aligned}
	\mathscr{F}''(t) &\geq -2I(u) > 4(d - J(u)) \geq 4(d - E(0)) =: C_4 > 0. 
\end{aligned}
\]
Since \(\mathscr{F}''(t)\geq C_4>0\), integrating twice gives $\mathscr{F}(t)>\frac{C_4}{2}t^2+\mathscr{F}'(0)t+\mathscr{F}(0)$.
By contradiction, if \(\|u(t)\|_{W_0}^2 \leq L_1\) for some \(L_1>0\) and all \(t\geq0\), then
$\mathscr{F}(t)\leq \frac{1}{2}C(n,s)L_1\,t+\|u_0\|_2^2$, which contradicts the quadratic growth of \(\mathscr{F}(t)\). Thus, \(\lim_{t\to\infty}\|u(t)\|_{W_0}^2=\infty\).

Case (II): \(E(0) = d\). We claim that there exists \(t_0 \in (0, T)\) such that
$
\int_0^{t_0} \|u_\tau\|_{W_0}^2 \mathrm{d}\tau > 0.
$
Suppose not, i.e., for all \(t \in [0, T)\),
$
\int_0^t \|u_\tau\|_{W_0}^2 \mathrm{d}\tau \equiv 0.
$
This implies \(\|u_\tau\|_{W_0} = 0\) for a.e. \(\tau\in[0,T)\). Since \(u_t(x) = 0\) for \(x\in \mathbb{R}^n\setminus \Omega\) and \(u_t \in C\left([0,T);L^2(\Omega)\right)\), it follows that \(u_t \equiv 0\) on \([0,T)\). Thus $u(t) \equiv u_0$ on $[0, T)$, so $I(u_0) = 0$, contradicting $I(u_0) < 0$. By \eqref{2.6}, 
\[
E(t_0) = d - \frac{1}{2}C(n,s)\int_0^{t_0} \|u_\tau\|_{W_0}^2 \mathrm{d}\tau < d.
\]
If we take $t_0$ as the new initial time and repeat the argument in Case (I), we obtain the same conclusion.
\hfill$\square$

{\bfseries Proof of \cref{t2.4}}

Proof of (i). Step 1: Extension of \(u\) to infinite time.

Since $(u_0, u_1) \in \mathcal{U}_{\mathscr{O}}$, from \cref{l3.7} we know $I(u(t)) < \|u_t(t)\|_2^2 - \mathscr{O}$ for \( t \in[0, T) \) and
\begin{eqnarray} \label{4.17}
	d < E(0) \leq \frac{\sqrt{\|(-\Delta)^{\frac{s}{2}}u(t)\|_{2}^2 + 2(u(t), u_t(t))}}{2\sqrt{e b_*}} - \frac{\mathscr{O}}{2}, \quad \forall t \in [0, T).
\end{eqnarray}
By \cref{l3.4}, $d > \frac{1}{4} \exp\left( \frac{n}{s}(1+\ln a_0) \right)$, $a_0 = \sqrt{\frac{2}{K(n,s)}}>0$. Since $K(n,s) < e^2 < 2e^2$, it follows that $d > \frac14$. Since $b_* \geq 2$, we have
\begin{eqnarray} \label{4.18}
	\|(-\Delta)^{\frac{s}{2}}u(t)\|_{2}^2 + 2(u(t), u_t(t)) \geq 8 e d^2 > \frac12 e >  1, \quad \forall t \in [0, T).
\end{eqnarray}
From \eqref{4.17} and \eqref{4.18}, one obtains 
\begin{eqnarray}  \label{4.19}
	d<E(0) \leq \frac{\left\|(-\Delta)^{\frac{s}{2}} u\right\|_{2}^2+2\left(u, u_{t}\right)}{2 \sqrt{e b_*}}, \quad t \in[0, T).
\end{eqnarray}

Define $H(t) = \|u(t)\|_2^2$. From \eqref{2.9},  Cauchy-Schwarz inequality and Young's inequality, we deduce
\begin{eqnarray} \label{4.20}
	H'(t) = 2(u(t), u_t(t)) \leq \|u(t)\|_2^2 + \left( 2 - \frac{2}{\sqrt{e b_*}} \right) \|u_t(t)\|_2^2.
\end{eqnarray}
Write $
N(s,n) = \sqrt{\frac{2}{e}}
-\frac{n}{s}\left[
1+\frac12\ln\left(\frac{2-\sqrt{2/e}}{K(n,s)}\right)
\right] + 1.
$ By \eqref{4.20}, \eqref{2.3}, \eqref{2.5}, \eqref{2.6}, \eqref{4.19}, \eqref{3.1}, the Cauchy-Schwarz inequality, Young's inequality and $b_* \geq 2$ from \eqref{2.9}, for any $a > 0$,
\[
\begin{aligned}
	&(e^{N(s,n)}+ 1 + 2H(t))\ln(e^{N(s,n)} + 1 + 2H(t)) - H'(t) \\
	& = H(t)\ln(e^{N(s,n)} + 1 + 2H(t)) + (e^{N(s,n)}+ 1 + H(t))\ln(e^{N(s,n)} + 1 + 2H(t)) - H'(t) \\
	& \geq N(s,n) \|u\|_2^2 + \|u\|_2^2\ln \|u\|_2^2 -2 \|u_t\|_2^2 + \frac{2}{\sqrt{e b_*}}\|u_t\|_2^2 - \|u\|_2^2 \\
	& \geq N(s,n) \|u\|_2^2 + \|u\|_2^2\ln \|u\|_2^2 + 4\left( \frac{1}{4}\|u\|_2^2 + \frac{1}{2}I(u) - E(0) \right) + \frac{2}{\sqrt{e b_*}}\|u_t\|_2^2 - \|u\|_2^2 \\
	& \geq N(s,n) \|u\|_2^2 + \|u\|_2^2\ln \|u\|_2^2 + 2\left\|(-\Delta)^{\frac{s}{2}} u\right\|_{2}^2 - 2\int_\Omega |u|^2 \ln|u|\,dx - \frac{2\left\|(-\Delta)^{\frac{s}{2}} u\right\|_{2}^2+4\left(u, u_{t}\right)}{\sqrt{e b_*}} + \frac{2}{\sqrt{e b_*}}\|u_t\|_2^2 \\
	& \geq N(s,n) \|u\|_2^2 +  2\left\|(-\Delta)^{\frac{s}{2}} u\right\|_{2}^2 + \frac{n}{s}(1+\ln a)\|u\|_2^2 - a^2 K(n,s) \left\|(-\Delta)^{\frac{s}{2}} u\right\|_{2}^2 - \frac{2\left\|(-\Delta)^{\frac{s}{2}} u\right\|_{2}^2}{\sqrt{e b_*}} - \frac{2}{\sqrt{e b_*}}\|u\|_2^2  \\
	& \geq N(s,n) \|u\|_2^2 + \left(2 - a^2 K(n,s) - \frac{2}{\sqrt{e b_*}}\right)\left\|(-\Delta)^{\frac{s}{2}} u\right\|_{2}^2 + \left(\frac{n(1+\ln a)}{s} - \frac{2}{\sqrt{e b_*}}\right)\|u\|_2^2 \\
	& \geq \left(2 - a^2 K(n,s) - \sqrt{\frac{2}{e}}\right)\left\|(-\Delta)^{\frac{s}{2}} u\right\|_{2}^2 + \left(N(s,n) + \frac{n(1+\ln a)}{s} - \sqrt{\frac{2}{e}}\right)\|u\|_2^2 \\
	& := C(a)\left\|(-\Delta)^{\frac{s}{2}} u\right\|_{2}^2 + D(a)\|u\|_2^2.
\end{aligned}
\]
Taking $a \in \left(
\sqrt{\frac{2-\sqrt{2/e}}{K(n,s)}}\,e^{-s/n},
\sqrt{\frac{2 - \sqrt{2/e}}{K(n,s)}}
\right)$, we find that $C(a) > 0$ and $D(a) > 0$. Consequently, we obtain $(1 + H(t))\ln(1 + H(t)) \geq H'(t)$. By Gronwall's inequality,
\begin{eqnarray} \label{4.21}
	\|u\|_2^2 = H(t) \leq \frac12\left[
	\big(e^{N(s,n)} + 1 + 2H(0)\big)^{e^{2t}}
	-\big(e^{N(s,n)} + 1\big)
	\right] =: H_2(t).
\end{eqnarray}
By an argument similar to that for \eqref{3.9}, combined with \eqref{4.21}, we obtain
\begin{eqnarray} \label{4.22}
	\begin{aligned}
		|E(0)| &\geq E(t) \geq \frac12\big\|u_t\big\|_2^2 + A(a) \left\|(-\Delta)^{\frac{s}{2}} u\right\|_{2}^2 + B(a) \big\|u\big\|_2^2 - \frac{1}{4} \|u\|_2^2 \ln\|u\|_2^2 \\
		& \geq \frac12\big\|u_t\big\|_2^2 + A(a) \left\|(-\Delta)^{\frac{s}{2}} u\right\|_{2}^2 + \frac{n(1+\ln a)}{4s} \big\|u\big\|_2^2 - \frac{1}{4}(H_2(t) + 1)\ln(1 + H_2(t)),
	\end{aligned}
\end{eqnarray}
where $A(a) = \frac{2 - a^2 K(n,s)}{4}$. Substituting $a = e^{-1}$ into \eqref{4.22}, we obtain $A(e^{-1}) > \frac{1}{4}$ since $K(n,s) < e^2$, and
\[
\left\|(-\Delta)^{\frac{s}{2}} u(t)\right\|_{2}^2 + \|u_t(t)\|_2^2 \leq \frac{1}{4A(e^{-1})} \left((H_2(t) + 1)\ln(1 + H_2(t)) + 4|E(0)|\right) =: \varpi_1(t).
\]
Since $\|u(t)\|_{W_0}^2 = \dfrac{2}{C(n,s)}\left\|(-\Delta)^{\frac{s}{2}}u(t)\right\|_{2}^2$, the above bound implies $\limsup_{t\to T^-}(\|u_t(t)\|_2^2 + \|u(t)\|_{W_0}^2) < \infty$. By \cref{t2.1}, we conclude \( T=\infty \).

Step 2: Polynomial growth. 

By \eqref{2.8} and \cref{l3.7}, $I(u) < \|u_t\|_2^2 - \mathscr{O}$. This, together with \eqref{4.16}, yields
$
\mathscr{F}''(t) = 2\|u_t\|_2^2 - 2I(u) > 2\mathscr{O},  t \geq 0.
$
Let
$
\hat{\mathscr{O}} = \min\left\{ \|u_1\|_2^2 - I(u_0), \frac{\sqrt{\|(-\Delta)^{\frac{s}{2}}u_0\|_{2}^2 + 2(u_0,u_1)}}{\sqrt{e b_*}} - 2E(0) \right\}.
$
By \eqref{2.8}, $(u_0, u_1) \in \mathcal{U}_{\mathscr{O}}$ for all $\mathscr{O} \in (0, \hat{\mathscr{O}})$.
As \(\mathscr{O} \to \hat{\mathscr{O}}\), the same argument based on \eqref{4.16} yields
$
\lim_{t \to \infty} \|u(t)\|_{W_0}^2 = \infty.
$

Proof of (ii). Fix any $\hat{Q} > d$. Choose $v \in W_0 \setminus \{0\}$ and $\psi \in L^2(\Omega) \setminus \{0\}$ with $(v, \psi) = 0$. Define
\[
u_0^{\hat{Q}} = \hat{m} v(x), \quad u_1^{\hat{Q}} = \hat{m} v(x) + \hat{h} \psi(x), \quad (\hat{m},\hat{h}>0 \text{ to be determined}).
\]
By orthogonality,
$$
(u_0^{\hat{Q}}, u_1^{\hat{Q}}) = \hat{m}^2 \|v\|_2^2, \quad \|u_1^{\hat{Q}}\|_2^2 = \hat{m}^2 \|v\|_2^2 + \hat{h}^2 \|\psi\|_2^2.
$$
From \eqref{2.2} and \eqref{2.5}, the initial energy satisfies
$
E(0) = \frac{\hat{h}^2}{2}\|\psi\|_2^2 + H(\hat{m}),
$
where
\[
H(\hat{m}) = \frac{\hat{m}^2}{2}\left( \frac{3}{2}\|v\|_2^2 + \frac{1}{2}C(n,s)\|v\|_{W_0}^2 - \int_\Omega v^2 \ln|v| dx - \ln \hat{m} \cdot \|v\|_2^2 \right).
\]

We now verify \eqref{2.8} and $E(0) = \hat{Q}$ for suitable $\hat{m},\hat{h}$. By \eqref{2.3},
\[
I(u_0^{\hat{Q}}) - \|u_1^{\hat{Q}}\|_2^2
= \hat{m}^2 \left( \frac{1}{2}C(n,s)\|v\|_{W_0}^2 - \ln \hat{m} \cdot \|v\|_2^2 - \int_\Omega v^2 \ln|v|\,dx - \|v\|_2^2 \right),
\]
which is negative for sufficiently large $\hat{m}$, so that $I(u_0^{\hat{Q}}) < \|u_1^{\hat{Q}}\|_2^2$. As $\lim_{\hat{m}\to\infty} H(\hat{m}) = -\infty$, take $\hat{m}$ large enough such that $\hat{Q} > H(\hat{m})$ and
$
d < \hat{Q} < \frac{\hat{m} \sqrt{\frac{1}{2}C(n,s)\|v\|_{W_0}^2 + 2\|v\|_2^2}}{2\sqrt{e b_*}}.
$
Setting
$
\hat{h} = \frac{\sqrt{2(\hat{Q} - H(\hat{m}))}}{\|\psi\|_2}
$
yields $E(0; (u_0^{\hat{Q}}, u_1^{\hat{Q}})) = \hat{Q}$. Thus $(u_0^{\hat{Q}}, u_1^{\hat{Q}})$ meets the blow-up condition. This completes the proof.
\hfill$\square$

\section*{Acknowledgments}
The authors would like to express their thanks to the referees for valuable advice regarding previous version of this paper. The research was supported by National Natural Science Foundation of China (12461021).

\section*{Competing interests}
The authors declare that they have no competing interests.

\end{document}